[1994/12/01]% LaTeX date must December 1994 or later\
\documentclass[12pt]{article}
\title{Two direct solvers for a system of linear equations}
\author{Michael F. Zimmer}
\date{\today}	

\usepackage{amsmath,amsthm,url}

\usepackage[utf8]{inputenc} % allow utf-8 input
\usepackage[T1]{fontenc}    % use 8-bit T1 fonts

\begin{document}
\maketitle

\begin{abstract}

A system of linear equations is normally understood as a linear mapping between two vector spaces.
However, most direct solutions (e.g., QR, LU, ...) rely on the inelegant approach
of back-substitution: a significant departure from such a characterization.
In this paper, two new methods are introduced which respect the underlying vector space
throughout, whether it be the row or column space of the coefficient matrix.
Both solutions produce a generalized inverse as well as a projection operator for the null space.
They also have the unique feature of admitting an online solution, which allows a solution
to be computed as the data becomes available. \\

\end{abstract}

%Two new algorithms for computing a direct solution to a system of linear equations are presented. 
%A variation on the orthonormalization step in the usual QR method allows one to bypass subsequent elimination or matrix inversion steps.
%The algorithms are formulated with respect to row and column spaces, and produce a generalized inverse in each case.  
%They may also be formulated in an online manner, so that the solution may be computed as the data for the equations become available. 

\noindent \emph{AMS classification:} 15A06, 15A09\\

\noindent \emph{Key Words:}  system of linear equations; row space; column space; generalized inverse; minimum norm; least squares; direct solution; online solution

% @@@@@@@@@@@@@@@@@@@@@@@@@@@@@@@@@@@@@@@@@@
\section{Introduction}

The problem of solving a system of linear equations
is widespread across mathematics, science and engineering.
Taking there to be $m$ equations in $n$ unknowns, the equations will be written as $Ax=b$, where $A=(a_{ij})$ is the coefficient matrix, $x=(x_j)$ is the vector of unknowns, and $b=(b_i)$.  
The indices are $i=1, \dots, m$ and $j=1, \dots ,n$, and the
entries in $A$, $x$ and $b$ belong to a field ${\cal F}$.
It will be assumed that ${\cal F}$ is endowed with an inner product;
here ${\cal F}$ is taken to be ${\cal C}$.

When the nullity of $A$ is non-zero, the inverse to $A$ does not exist, and it will no longer be possible to write the solution in the cogent form $x = A^{-1}b$. 
Nevertheless, in this case one may still use the generalized inverse ($G$) to express the particular and homogeneous parts of the solution as: 
\begin{align*}
x_p & = Gb \\
x_h & = Py \\
P & = 1_n - GA
\end{align*}
where $Ax_p=b$, $Ax_h=0$, $1_n$ is an n-by-n identity matrix, $P$ is a null space projection operator, and $y \in {\cal F}^n$.  
By definition, $Py \in {\cal N}[A]$, where ${\cal N}[A]$ is the null space of $A$.

The methods to be discussed will be seen to bear a similarity to the family of $QR$ factorization methods for solving $Ax=b$ (i.e., the $QR$, $RQ$, $LQ$ and $QL$ methods)\cite{watkins}.  
Those methods involve, for example, making the substitution $A=QR$ in $Ax=b$, and then solving for $x$.  
In contrast, the methods proposed herein involve a factorization of $A$ while it is in the context of the equation $Ax=b$.  
They'll be referred to as {\em in situ} methods, to emphasize this point.
In contrast, since the $QR$ methods involve a factorization of $A$ done externally to the context of $Ax=b$, 
they might instead be called {\em ex situ} methods.
In addition, it will be shown that in situ approaches allow for an {\em online} solution, which is not prima facie possible with the $QR$ methods.
Finally, two appendices are included: one for computing the Penrose identities for each method, and one to
comment on an {\em in situ} version of LU decomposition.

% @@@@@@@@@@@@@@@@@@@@@@@@@@@@@@@@@@@@@@@@@@
\section{Definitions}

The two main algorithms to be presented will involve orthonormalizing the row (column) vectors of the coefficient matrix $A$, by left (right) multiplying $A$ with a series of nonsingular matrices $M_s$ ($s=1,2,\dots $).  It's not necessary for the $M_s$ to implement elementary row (column) operations, but it's convenient to view them that way.  In the course of applying an orthonormalization procedure (e.g., Gram-Schmidt \cite{golub,langou,watkins,CGSMGS}) to the rows or columns of $A$, it may be that a zero vector will appear, since there may be linear dependencies amongst those vectors. If they do appear, they will be left in place in $A$ and thereafter ignored.  The result will be that the vectors of the transformed $A$ will form a quasi-orthonormal list of vectors, as defined next.\\

\noindent{\it Definition}: A {\bf quasi-orthonormal} list of vectors consists of vectors with (Euclidean) norm equal to $1$ or $0$.
Those vectors of norm $1$ are mutually orthogonal.\\ 

It will later be convenient to use a matrix which is essentially a p-by-p identity matrix, but with some of the $1$s replaced by $0$s.  
The row (or column) indices of the $1$s will comprise the set $S$.  Note that $S \subseteq \{ 1,2, \dots , p \}$.
Thus an {\bf index matrix} is defined as:
\[({\cal I}_p^S)_{ij}= \begin{cases}
1 & \text{if $i=j$ and $i \in S$}, \\
0 & \text{otherwise.}
\end{cases} \]
where $i,j$ each range from $1$ to $p$.
Note that if $|S| = p$, then ${\cal I}_p^S$ is just $1_p$.

% @@@@@@@@@@@@@@@@@@@@@@@@@@@@@@@@@@@@@@@@@@
\section{Method: {\em row space in situ}}
\label{sec:row}

The first step to the {\em row space in situ} method is to attempt to orthonormalize the rows of $A$, using the approach outlined in the previous section.
The result is that $Ax=b$ becomes
\begin{align*}
(\cdots M_2 M_1) A x & = (\cdots M_2 M_1)b \\
A'x & = b'
\end{align*}
where $M = (\cdots M_2 M_1)$, $A'=MA$, and $b'=Mb$.
The row vectors are now a quasi-orthonormal list of vectors.

Next define the set $S$ to contain the indices of the rows of $A'$ that are non-zero. 
It is then easy to verify
\begin{align}\label{AAstar}
A'~ (A')^*   & =  {\cal I}_m^S \\
{\cal I}_m^S ~ A'  & =  A' \\
{\cal I}_m^S ~ b'  & =  b'
\end{align}
where the superscript * represents a conjugate transpose.
A solution will now be claimed and verified. Afterwards, a constructive approach will be taken.\\

\noindent{\it Theorem}:
A solution to $Ax=b$ is 
\begin{equation}\label{mainsoln}
x = (A')^*b' + x_h
\end{equation}
where the variable definitions from this section are used.\\

\begin{proof}
The proposed solution is verified by substitution into $Ax$
\begin{align*}
Ax & = A(A')^* b' + Ax_h \\
   & = M^{-1} A'(A')^* b' \\
   & = M^{-1} {\cal I}_m^S b' \\
   & = M^{-1} b' \\
   & = b
\end{align*}
\end{proof}

Having verified an unmotivated solution to $Ax=b$, the solution will now be derived by constructive means.\\

\noindent{\it Lemma}:
An arbitrary vector $x \in {\cal F}^n$ may be expressed as
\begin{equation}\label{paramsoln}
x = (A')^*w + x_h
\end{equation}
where $A$ is an m-by-n matrix, $w \in {\cal F}^m$ and $x_h$ satisfies $Ax_h = 0$.\\

\begin{proof}
Because of the linear transformation $A': {\cal F}^n \rightarrow {\cal F}^m$, it is known (Thm. 18.3 \cite{bowen}) that
\begin{equation*}
{\cal F}^n = {\cal R}[ (A')^*] \oplus {\cal N}[A']
\end{equation*}
where ${\cal R}[A]$ is the range of $A$.
This direct sum decomposition may be used to rewrite an $x \in {\cal F}^n$ as
\begin{equation*}
x = w_1 + w_2
\end{equation*}
where 
\begin{align*}
w_1 & \in {\cal R}[ (A')^*] \\
w_2 & \in {\cal N}[A']
\end{align*}
However, every vector that is in the range of some matrix $Q$ may be expressed as $Qv$,
for an appropriate $v$.  Thus, for some $w \in {\cal F}^m$, $w_1$
may be written as 
\begin{equation*}
w_1 = (A')^* w
\end{equation*}
Also, because $M$ is nonsingular, ${\cal N}[A'] = {\cal N}[A]$,
so that $w_2 \in {\cal N}[A]$.
\end{proof}

What remains is to determine $w$ and to find a means of computing $x_h$.
Note that upon substituting \eqref{paramsoln} into $A'x=b'$ 
one obtains $A'(A')^*w = b'$, 
which upon using \eqref{AAstar} becomes ${\cal I}_m^S ~w = b'$. 
This can be used to re-express equation \eqref{paramsoln} as 
\begin{align*}
x & = (A')^* ~w + x_h\\
  & = (A')^* {\cal I}_m^S ~w + x_h\\
  & = (A')^* ~b' + x_h
\end{align*}
which is the same as in equation \eqref{mainsoln}.

% =============================================================
\subsection*{Discussion}

Of note is that $x_p = Gb = (A')^* b'$ is a minimum norm solution.
This follows since it is formed from row space vectors, and not null space vectors.  
Its minimum norm nature is also evident since the generalized inverse is, in general, of type $\{124\}$ (see Appendix).

There are two variations to this approach, which are mainly of interest for numerical applications.
When $G$ and $M$ are not explicitly computed, the solution is:
\begin{align*}
x_p & = (A')^* b' \\
x_h & = Py \\
P & = 1_n - (A')^* A'
\end{align*}
where again $y$ is an arbitrary vector in ${\cal F}^n$.
On the other hand, if $G$ is explicitly computed, the solution may be written as:
\begin{align*}
x_p & = Gb \\
x_h & = Py \\
G & = (A')^* M \\
P & = 1_n - GA
\end{align*}
where $y$ is the same.  
Depending on the application, it may or may not be advantageous to explicitly compute $G$.

% @@@@@@@@@@@@@@@@@@@@@@@@@@@@@@@@@@@@@@@@@@
\section{Method: {\em column space in situ}}
\label{sec:col}

In this {\em column space in situ} method, it's first recognized that $Ax=b$ will have a solution only when $b$ lies in the column space of $A$.  
In particular, a least squares solution is sought, in which only the portion of $b$ that lies in that column space is allowed; this is defined as $b_c$.  
Hence, the equation to be solved is $Ax=b_c$, and the strategy now is to solve it by attempting to orthonormalize the columns of $A$
by right-multiplying it by $M_s$.
Thus, matrix products $M_s M_s^{-1}$ are repeatedly inserted between $A$ and $x$ in $Ax=b_c$, to produce
\begin{equation*}
A (M_1 M_2 \cdots ) (\cdots M_2^{-1} M_1^{-1}) x = b_c
\end{equation*}
or
\begin{equation*}
A' M^{-1} x = b_c
\end{equation*}
where $M = (M_1 M_2 \cdots )$ and $A' = AM$.\\

As explained in the section Definitions, the matrices $M_s$ will effect an orthonormalization procedure on the columns of $A$.  Observe that $A'$ can be used in a projection operator to compute  
$b_c = A' (A')^* b$, and the equation to be solved can now be written
\begin{align*}
A' M^{-1} x & = A' (A')^* b \\
A [ x & - M(A')^*b ] = 0
\end{align*}
The piece in the square brackets must lie in ${\cal N}[A]$; it is named $x_h'$.  
Note that the solution $x$ and vector $x_h'$ are now entwined in the equation  $x - M(A')^*b = x_h'$ (and one may be used to define the other, depending on which is taken as given).  
Finally, observe that $x_h'$ may be defined to be any null space vector, and that it will simply shift the solution $x$.  
A solution is now proved which assumes an arbitrary null space vector $x_h$.\\

\noindent{\it Theorem}:
A solution to $Ax=b_c$ is 
\begin{equation}
x = M(A')^* b + x_h
\end{equation}
where the variable definitions from this section are used.
\begin{proof}
The proposed solution is verified by substitution into $Ax$
\begin{align*}
Ax & = AM(A')^* b + Ax_h\\
   & = A'(A')^* b \\
   & = b_c
\end{align*}
\end{proof}

% =============================================================
\subsection*{Discussion}

In this column space approach, the set $S$ is defined to be the columns of $A'$ that are nonzero, so that ${\cal I}_n^S$ is the index matrix associated with $A'$.  It is easy to verify:
\begin{align}\label{AAstar_ids_col}
A'		& = A' {\cal I}_n^S  \\
(A')^*  & = {\cal I}_n^S (A')^*
\end{align}
\\
These identities are useful in determining that the generalized inverse $G = M(A')^*$ is of type $\{123\}$ (see Appendix).  In summary, the (least squares) solution found for this section can rewritten as
\begin{align*}
x_p & = Gb \\
x_h & = Py \\
P  & = 1_n - GA
\end{align*}
where once again, $y$ is an arbitrary vector in $F^n$.

% @@@@@@@@@@@@@@@@@@@@@@@@@@@@@@@@@@@@@@@@@@
\section{Online capabilities}

It is normally assumed that all data (i.e., $A$,$b$) is present before solving $Ax=b$ may commence, using any of the usual methods.
For the algorithms introduced here, it will be seen this is {\em not} required.
Indeed, the linear solver may begin its work when the data arrives in a certain regular manner;
this is referred to as an {\em online} algorithm.
This capability may find use in the case where the dimensions of $A$ are large, and it is also a time-consuming process to compute its entries.

% =============================================================
\subsection*{Row space}

Assume that the rows of $[A|b]$ become available one at a time to the solver, in the order $1$ to $m$.  
The orthonormalization steps described earlier may be done incrementally as these rows become available.  
When the GS procedure is used to effect the orthonormalization, the $1$st through $i$th rows of $A$, $M$ and $b$ will no longer change following the $i$th step of the algorithm.  
It follows that these rows can be used in a computation of $x_p$.
Taking advantage of the column-row expansion \cite{carlson}, 
$x_p$ may be written as
\begin{equation*}
x_p = (A')^* b' = \sum_{i=1}^m x_p^{(i)}
\end{equation*}
where
\begin{equation*}
x_p^{(i)} = \text{Col}_i[(A')^*] ~ b_i'
\end{equation*}
and Col$_i$ signifies the $i$th column.
Thus following the $i$th step in the algorithm, $x_p^{(i)}$ may be computed.  
Also note that the updates to $x_p$ are mutually orthogonal, i.e.,
$<x_p^{(i)}, x_p^{(k)} > = 0$ for $i \neq k$.
This means that the estimation of $\|x\|$ is
monotonically non-decreasing as more rows are included.

One can take a similar approach when the second variation is used.  For example, the generalized inverse may be computed in an online mode as:
\begin{equation*}
G = (A')^* M = \sum_{i=1}^m \text{Col}_i[(A')^*] ~ \text{Row}_i[M]
\end{equation*}
Note that the summand may be computed following the $i$th step of the GS procedure. 
Following the computation of $G$, the particular solution is easily computed.  \\

% =============================================================
\subsection*{Column space}

Here it will be assumed that all of $b$ is available at the outset, and that the columns of $A$ become available one at a time to the solver, in the order $1$ to $n$.  
When the GS procedure is used to effect the orthonormalization, the $1$st through $j$th columns of $A$ and $M$ will no longer change following the $j$th step of the algorithm.  
Since these columns in $M$ and $A$ are done changing at these points, they become available to be used in a computation of $x_p$. 

Again, a column-row expansion \cite{carlson} may be used to rewrite $x_p$, this time as
\begin{align*}
x_p & = \sum_{j=1}^n x_p^{(j)} \\
x_p^{(j)} & = \text{Col}_j[M] ~ \text{Row}_j[(A')^*] ~ b 
\end{align*}
Since the $j$th row of $(A')^*$ and $j$th column of $M$ are available at the $j$th step, $x_p^{(j)}$ can be computed at the $j$th step. Also observe that this same online technique can be trivially extended to compute $G$, $P$ or even $x_h$ in an online manner.

% =============================================================
\subsection*{Complexity of online approaches}

At play are two time scales: the time ($\tau$) needed for a new portion of the data to be acquired, and the time ($\Delta$) needed for the linear solver to process the new portion of data. 
The {\em additional (time) complexity} due to the solver is defined as the additional computations that must yet be performed after all the data has arrived.
As the orthonormalization steps are performed by the linear solver,
the complexity it incurs increases from $\Delta_{min}$ for the first step, to $\Delta_{max}$ for the final step.
If $\tau > \Delta_{max}$, the solver will generally be waiting on data acquisition.
If $\tau < \Delta_{min}$, the complexity will in general be comparable to when all data is available at the outset.
If $\tau =0$, the complexity will be the same as if all data is available at the outset.

For the row space approach, the time used by the solver on the $i$th row is $\Delta(i) \sim O(4ni-n)$, which varies from $\Delta_{min} = \Delta(1)$ to $\Delta_{max} = \Delta(m)$.  
Also, no matter how large the size of $\tau$, the solver will always add complexity due to the computations done on the last row ($i=m$), as well as the final matrix-vector multiply (i.e., forming $(A’)^* b’$).  
Combined, these two contribute a complexity of $O(mn)$, which acts as the minimum amount of added complexity for the row space solver.
The maximum complexity occurs when $\tau \ll \Delta_{min}$, and in this case is $O(m^2n)$.  
Thus the complexity added by the solver is:
\begin{equation*}
%O(4mn) < \text{complexity} < O(2m^2n)
O(mn) \leq \text{complexity} \leq O(m^2n)
\end{equation*}

Similarly, in the column space approach, the time used by the solver on the $j$th column is $\Delta(j) \sim O(4mj-m)$, which varies from $\Delta_{min} = \Delta(1)$ to $\Delta_{max} = \Delta(n)$.  
By a similar set of arguments, the complexity added by the solver in this case is
\begin{equation*}
%O(4mn) < \text{complexity} < O(2mn^2)
O(mn) \leq \text{complexity} \leq O(mn^2)
\end{equation*}
Thus the acquisition time $\tau$ acts as a crossover parameter for the complexity of the algorithm.

% @@@@@@@@@@@@@@@@@@@@@@@@@@@@@@@@@@@@@@@@@@
\section{Final remarks}

The approaches discussed herein, which came about while working on the inverse kinematics problem in robotics, takes a different view on factorization methods.  
Here, the factorization is done when $A$ is in the context of the equation $Ax=b$.  
Borrowing a term used in metallurgy, it is described as an {\em in situ} approach.  

One of the new features the algorithms present is that it can be computed in an online manner.  
This is expected to be useful in applications where the data is large, or when it simply takes a long time to create and/or acquire the data.
What is also immediately noteworthy about the algorithm is that it doesn't use an elimination or partitioning strategy; in particular, there was no solving of a triangular system of equations.

The two in situ approaches have a symmetry between them.  
The (second variation of the) {\em row space in situ} method involves a transformation of $[A|1]$ to produce a particular solution.
\begin{equation*}
[ A \vert 1 ] \longrightarrow [A' \vert M ] \longrightarrow x_p=(A')^* M \, b
\end{equation*}
Likewise, the {\em column space in situ} method may be expressed as:
\begin{equation*}
\begin{bmatrix}
 A \\
   1  
\end{bmatrix}
\longrightarrow
\begin{bmatrix}
A' \\
 M
\end{bmatrix}
\longrightarrow
x_p=M(A')^* \, b
\end{equation*}
The main reason for mentioning it is aesthetics: the two solutions comprise a more symmetrical viewpoint on solving $Ax=b$.

Although the new algorithms were cast to solve $Ax=b$, they can
easily solve its matrix generalization \cite{private_notes}: $AX=B$, where $X$ is n-by-p and $B$ is m-by-p (and $p \geq 1$).  
In that case the orthonormalization steps proceed as before,
and the particular and homogeneous parts of the solution 
$X = X_p + X_h$ are
\begin{align*}
 X_p & = G B \\
 X_h & = P Y 
\end{align*}
where $G$ and $P$ are the same as before, and $Y \in {\cal F}^{n \times p}$.
This also admits an online formulation.

% @@@@@@@@@@@@@@@@@@@@@@@@@@@@@@@@@@@@@@@@@@
\section*{Acknowledgements}
The author dedicates this paper to Robert E. Zimmer for his 
kind support and encouragement.  In addition, the author thanks 
Prof. Daniel Grayson for a discussion of the material.

% @@@@@@@@@@@@@@@@@@@@@@@@@@@@@@@@@@@@@@@@@@
% @@@@@@@@@@@@@@@@@@@@@@@@@@@@@@@@@@@@@@@@@@
\section*{APPENDIX: Penrose conditions}

It is convenient to classify a generalized inverse ($G$) according to which Penrose conditions it satisfies \cite{penrose1,penrose2,campbell,benis,bapat}. 
In particular, different properties of the solution $Gb$ follow if certain sets of the following Penrose conditions are satisfied. 

\begin{align*}
\text{(1)} \, & AGA = A \\
\text{(2)} \, & GAG = G \\
\text{(3)} \, & AG = (AG)^* \\
\text{(4)} \, & GA = (GA)^* 
\end{align*}

% =============================================================
\subsection*{Row space in situ}

Here the index matrix will be ${\cal I}_m^S$, where the set $S$ is the set of non-zero row vectors in the matrix $A'$, as defined in Sec.~\ref{sec:row}.  
The identities involving ${\cal I}_m^S$ are used herein.

The first Penrose condition is seen to always be true for our solution.
\begin{align*}
AGA & = A(A')^*MA  = M^{-1}MA(A')^*A'  \\
    & = M^{-1}A'(A')^*A' = M^{-1} {\cal I}_m^S A'  \\
    & = M^{-1}A' = A 
\end{align*}
Likewise, the second condition is always true.
\begin{align*}
GAG & = [(A')^*M]A[(A')^*M] \\
    & = (A')^*A'(A')^*M  = (A')^* {\cal I}_m^S M \\
    & = (A')^*M  = G
\end{align*}
The third Penrose condition is seen to be problematic
\begin{align*}
AG & = A(A')^* M  = M^{-1}MA(A')^*M \\
   & = M^{-1}A' (A')^*M  = M^{-1} {\cal I}_m^S M
\end{align*}
If $A$ is of full row rank, then ${\cal I}_m^S$ equals $1_m$,  $AG$ becomes the identity matrix $1_m$, and the condition is manifestly satisfied.  
However, if $A$ is \emph{not} of full row rank, then this condition is not in general satisfied.  
Finally, the fourth Penrose condition is seen to be true:
\begin{align*}
GA & = (A')^* MA = (A')^*A' \\
   & = ((A')^* A')^* = ((A')^* MA)^* \\
   & = (GA)^*
\end{align*}
In summary, the generalized inverse obtained by this algorithm is at least a $\{124\}$-inverse.  
If it's additionally true that $A$ is of full row rank, then it becomes a unique $\{1234\}$-inverse.
Recall that generalized inverses which are at least
of type $\{14\}$ (which is the case here) provide
solutions $Gb$ with a minimum-norm property.

% =============================================================
\subsection*{Column space in situ}

Here the index matrix will be ${\cal I}_n^S$, where the set $S$ is the set of non-zero column vectors in the matrix $A'$, as defined in Sec.~\ref{sec:col}.
The identities involving ${\cal I}_n^S$ are used herein.

The first Penrose condition is seen to always be true.
\begin{align*}
AGA & = A M (A')^* A = A' (A')^* A = A 
\end{align*}
To understand this, recall that $A'(A')^*$ projects away vectors
not in ${\cal R}(A)$; hence it leaves $A$ unchanged.
Likewise, the second condition is always true.
\begin{align*}
GAG & = [M (A')^*] A [M (A')^*] \\
    & = M (A')^* A' (A')^* = M {\cal I}_n^S (A')^* \\
    & = M  (A')^*  = G
\end{align*}
The third Penrose condition is also true
\begin{align*}
AG & = A M (A')^*  = A' (A')^* \\
   & = [ A' (A')^* ]^* = [ AM (A')^* ]^* \\
   & = ( AG )^* 
\end{align*}
The fourth Penrose condition $GA = (GA)^*$ is seen to
be problematic
\begin{align*}
GA & = M(A')^* A  = M (A')^* A MM^{-1} \\
   & = M (A')^* A' M^{-1}  = M {\cal I}_n^S M^{-1}
\end{align*}
Observe that if $A$ is of full column rank, then ${\cal I}_n^S$ equals the identity matrix $1_n$, and $GA$ becomes $1_n$; 
this trivially allows the identity to be satisfied.  
However, if $A$ is \emph{not} of full column rank, then this condition is not true in general.

In summary, the generalized inverse obtained by this
algorithm is at least a $\{123\}$-inverse.  If it's 
additionally true that $A$ has full column rank, then 
it becomes the unique $\{1234\}$-inverse.
Recall that generalized inverses which are at least
of type $\{13\}$ (which is the case here) provide
solutions $Gb$ with a least-squares property.

% @@@@@@@@@@@@@@@@@@@@@@@@@@@@@@@@@@@@@@@@@@
% @@@@@@@@@@@@@@@@@@@@@@@@@@@@@@@@@@@@@@@@@@
\section*{APPENDIX: LU decomposition}

The approach of LU decomposition can be reworked using the ideas of an in situ
factorization.  First, recall that this approach normally involves determining a lower (L) and upper (U)
triangular matrix factorization of the coefficient matrix as $A=LU$.  
This factorization is done externally (i.e., {\em ex situ}) to the context of the equation $Ax=b$.
This is then solved as two triangular sets of equations: first $Ly = b$, then $Ux=y$.  

The {\em in situ} version of LU decomposition is to again modify $A$ through left-multiplication through suitable matrices.
In this case the matrices $M_s$ cause $A'$ to have a triangular form.
%(Also, pivoting is easily implemented, for this or any of the methods in this paper.)
Since this is done in situ, these same operations are applied to the right-hand side of $Ax=b$.
While these operations can be applied directly to $b$ (i.e., as Gaussian elimination), 
it's preferred to equivalently apply them to $1_m b$ for more versatile use.  
Thus, $Ax= 1_m b$ becomes $A'x=Mb$, where again $M$ represents their cumulative product
$M = (\cdots M_2 M_1)$.  Once this is accomplished, there remain
two main steps in solving for $x$: computing $Mb$; solving $A'x=Mb$ using back substitution.
Note that $A$ is not required to be square or have full rank.

Finally, in some cases the equation $A'x=Mb$ may be interpreted as $Ux = L^{-1}b$,
in which case it might be called an {\em in situ LU decomposition}.

% @@@@@@@@@@@@@@@@@@@@@@@@@@@@@@@@@@@@@@@@@@
% @@@@@@@@@@@@@@@@@@@@@@@@@@@@@@@@@@@@@@@@@@


\begin{thebibliography}{10}
\bibitem{CGSMGS}
Note that the GS procedure can be implemented as either Classical (CGS) or Modified (MGS).  An analysis of their stability properties while using reorthogonalization is done in \cite{giraud}.  Also, see \cite{langou} for a review of Laplace's orthogonalization method, which appeared earlier.
\end{thebibliography}
\end{document}